\documentclass{amsart}
\usepackage{amsmath,amssymb}
\usepackage{yhmath}

\newcommand{\bdis}{\begin{displaymath}}
\newcommand{\edis}{\end{displaymath}}
\newcommand{\be}{\begin{equation}}
\newcommand{\ee}{\end{equation}}
\newcommand{\mbb}{\mathbb}
\newcommand{\mcal}{\mathcal}

\newcommand{\vp}{\varphi}

\newcommand{\zf}{\zeta\left(\frac{1}{2}+it\right)}

\DeclareMathOperator{\re}{Re}

\theoremstyle{definition}

\theoremstyle{remark}
\newtheorem{remark}[]{Remark}

\newtheorem*{mydef11}{{\bf Theorem 1}}

\newtheorem*{mydef12}{{\bf Theorem 2}}

\newtheorem*{mydef13}{{\bf Theorem 3}}

\newtheorem*{mydef51}{{\bf Lemma 1}}

\newtheorem*{mydef52}{{\bf Lemma 2}}

\newtheorem*{mydef53}{{\bf Lemma 3}}

\newtheorem*{mydef54}{{\bf Lemma 4}}

\newtheorem*{mydef55}{{\bf Lemma 5}}

\newtheorem*{mydef56}{{\bf Lemma 6}}

\newtheorem*{mydef57}{{\bf Lemma 7}}

\newtheorem*{mydef81}{{\bf Property 1}}

\newtheorem*{mydef91}{{\bf Formula1}}

\newtheorem*{mydef92}{{\bf Formula2}}

\newtheorem*{mydef93}{{\bf Formula3}}

\numberwithin{equation}{section}

\begin{document}

\title[Jacob's ladders, almost linear increments \dots ]{Jacob's ladders, almost linear increments of the Hardy-Littlewood integral (1918), the classical Dirichet's sum of divisors (1849) and their relationship with the Fermat-Wiles theorem}

\author{Jan Moser}

\address{Department of Mathematical Analysis and Numerical Mathematics, Comenius University, Mlynska Dolina M105, 842 48 Bratislava, SLOVAKIA}

\email{jan.mozer@fmph.uniba.sk}


\begin{abstract}
In this paper we obtain number of new equivalents of the Fermat-Wiles theorem that are based on Jacob's ladders. The main of these is the $D$-equivalent that is generated by the Dirichlet's $D(x)$-function.  
\end{abstract}
\maketitle

\section{Introduction} 

\subsection{}    

Let us remind that we have obtained in our previous paper \cite{8} the point of contact between almost linear increments of the Hardy-Littlewood integral and the Fermat-Wiles theorem. Namely, we have proved, for example, the following: the $\zeta$-condition 
\be \label{1.1} 
\lim_{\tau\to\infty}\frac{1}{\tau}\int_{\frac{x^n+y^n}{z^n}\frac{\tau}{1-c}}^{[\frac{x^n+y^n}{z^n}\frac{\tau}{1-c}]^1}\left|\zf\right|^2{\rm d}t\not= 1
\ee  
on the class of all Fermat's rationals 
\be \label{1.2} 
\frac{x^n+y^n}{z^n},\ x,y,z,n\in\mbb{N},\ n\geq 3, 
\ee 
represents the $\zeta$-equivalent of the Fermat-Wiles theorem. In this paper we shall prove some new results in this direction. 

Next, let us remind also the Dirichlet's function $D(x)$ defined by the sum 
\be \label{1.3} 
D(x)=\sum_{n\leq x}d(n) 
\ee  
where $d(n)$ denotes the number of divisors of $n$ and 
\be \label{1.4} 
D(x)=D(N),\ x\in [N,N+1),\ N\in\mbb{N}. 
\ee 
Dirichlet himself have proved in an elementary way the following formula (Abhand. Al. Wiss., Berlin, 1849) 
\be \label{1.5} 
\begin{split}
& D(x)=x\ln x+(2c-1)x+\Delta(x), \\ 
& \Delta(x)=\mcal{O}(\sqrt{x}). 
\end{split}
\ee 

\subsection{} 

First, we have obtained the following result. 

The $D$-condition 
\be \label{1.6} 
\lim_{\tau\to\infty}\frac{1}{\tau}
\left\{
D\left(\left[\frac{x^n+y^n}{z^n}\frac{\tau}{1-c}\right]^1\right)-D\left(\frac{x^n+y^n}{z^n}\frac{\tau}{1-c}\right) 
\right\}\not=1, 
\ee 
where 
\be \label{1.7} 
\left[\frac{x^n+y^n}{z^n}\frac{\tau}{1-c}\right]^1=\vp_1^{-1}\left(\frac{x^n+y^n}{z^n}\frac{\tau}{1-c}\right), 
\ee  
on the class of all Fermat's rationals represents the $D$-equivalent of the Fermat-Wiles theorem. 

\begin{remark}
To gain the above mentioned result we have used the elementary Dirichlet's formula (\ref{1.5}) only. No improvements\footnote{For known improvements of the error term $\Delta(x)$ see for example \cite{9}, pp. 262 -- 275.} of the error term $\Delta(x)$ were used. 
\end{remark} 

\begin{remark}
The above mentioned $D$-equivalent can be viewed as a kind of quadruple point of contact between foursome objects: 
\begin{itemize}
	\item[(a)] The Riemann's zeta function on the critical line, 
	\item[(b)] Jacob's ladder $\vp_1(t)$\footnote{See \cite{2}.}, 
	\item[(c)] Dirichlet's function $D(x)$, 
	\item[(d)] Fermat-Wiles theorem. 
\end{itemize}
\end{remark} 

\subsection{} 

Second, we have obtained the following result. 

The $\zeta$-condition 
\be \label{1.8} 
\lim_{\rho\to\infty}\frac{1}{\rho}\int_{\vp_1(\frac{x^n+y^n}{z^n}\frac{\rho}{(1-c)\sigma(l)})}^{\frac{x^n+y^n}{z^n}\frac{\rho}{(1-c)\sigma(l)}}\left\{ \sigma(l)\left|\zf\right|^2+(1-c)|S_1(t)|^{2l}\right\}{\rm d}t\not=1, 
\ee  
where 
\be \label{1.9} 
S_1(t)=\frac{1}{\pi}\int_0^t \arg\zf{\rm d}t, 
\ee  
on the class of all Fermat's rationals represents the next $\zeta$-equivalent of the Fermat-Wiles theorem. 

\subsection{} 

Third, we have obtained the following result. 

The $D$-condition 
\be\label{1.10} 
\lim_{\tau\to\infty}\ln\left\{
D\left(\left[\exp\left(\frac{x^n+y^n}{z^n}\ln \tau\right)\right]^1\right)-
D\left(\exp\left(\frac{x^n+y^n}{z^n}\ln \tau\right)\right)
\right\}^{\frac{1}{\ln\tau}}\not=1
\ee 
on the class of all Fermat's rationals represents the next $D$-equivalent of the Fermat-Wiles theorem. 

\subsection{} 

In this paper we use the following notions from our papers \cite{2} -- \cite{5}: 
\begin{itemize}
	\item[{\tt (a)}] Jacob's ladder $\vp_1(t)$, 
	\item[{\tt (b)}] the function 
	\bdis 
	\begin{split}
	& \tilde{Z}^2(t)=\frac{{\rm d}\vp_1(t)}{{\rm d}t}=\frac{1}{\omega(t)}\left|\zf\right|^2,\\ 
	& \omega(t)=\left\{1+\mcal{O}\left(\frac{\ln\ln t}{\ln t}\right)\right\}\ln t,\ t\to\infty, 
	\end{split}
	\edis 
	\item[{\tt (c)}] direct iterations of Jacob's ladders 
	\bdis 
	\begin{split}
	& \vp_1^0(t)=t,\ \vp_1^1(t)=\vp_1(t),\ \vp_1^2(t)=\vp_1(\vp_1(t)),\dots , \\ 
	& \vp_1^k(t)=\vp_1(\vp_1^{k-1}(t))
	\end{split}
	\edis 
	for every fixed natural number $k$, 
	\item[{\tt (d)}] reverse iterations of Jacob's ladders 
	\bdis 
	\begin{split}
	& \vp_1^{-1}(T)=\overset{1}{T},\ \vp_1^{-2}(T)=\vp_1^{-1}(\overset{1}{T})=\overset{2}{T},\dots, \\ 
	& \vp_1^{-r}(T)=\vp_1^{-1}(\overset{r-1}{T})=\overset{r}{T},\ r=1,\dots,k, 
	\end{split} 
	\edis 
	that is, for example, 
	\be \label{1.11} 
	\vp_1(\overset{r}{T})=\overset{r-1}{T}
	\ee  
	for every fixed $k\in\mbb{N}$, where 
	\be \label{1.12}
	\begin{split} 
	& \overset{r}{T}-\overset{r-1}{T}\sim(1-c)\pi(\overset{r}{T});\ \pi(\overset{r}{T})\sim\frac{\overset{r}{T}}{\ln \overset{r}{T}},\ r=1,\dots,k,\ T\to\infty, \\ 
	& \overset{0}{T}=T<\overset{1}{T}(T)<\overset{2}{T}(T)<\dots<\overset{k}{T}(T), \\ 
	& T\sim \overset{1}{T}\sim \overset{2}{T}\sim \dots\sim \overset{k}{T},\ T\to\infty. 
	\end{split}
	\ee 
\end{itemize}  

\begin{remark}
The asymptotic behaviour of the points 
\bdis 
\{T,\overset{1}{T},\dots,\overset{k}{T}\}
\edis  
is as follows: at $T\to\infty$ these points recede unboundedly each from other and all together are receding to infinity. Hence, the set of these points behaves at $T\to\infty$ as one-dimensional Friedmann-Hubble expanding Universe. 
\end{remark} 

\section{On the point of contact between the Riemann's zeta-function on the critical line and the set of $L_2$-orthogonal systems}

Let us remind that our $D$-results (\ref{1.6}) and (\ref{1.10}) give new point of contact between completely different objects. Namely, between the Dirichlet's $D(x)$-function and the Fermat-Wiles theorem. 

For completeness, we remind also our result, see \cite{6}, that gives the first point of contact in this direction. Namely, the point of contact between the Riemann's zeta-function on the critical line and the theory of $L_2$-orthogonal systems. We shortly recapitulate this result in the following text. 

\subsection{} 

We have introduced the generating vector operator $\hat{G}$ acting on the class of all $L_2$-orthogonal systems 
\be \label{2.1} 
\{f_n(t)\}_{n=0}^\infty,\ t\in [a,a+2l],\ a\in\mbb{R},\ l>0
\ee 
as 
\be \label{2.2} 
\begin{split}
& \{ f_n(t)\}_{n=0}^\infty\xrightarrow{\hat{G}}\{ f^{p_1}_n(t)\}_{n=0}^\infty\xrightarrow{\hat{G}}\{ f^{p_1,p_2}_n(t)\}_{n=0}^\infty\xrightarrow{\hat{G}} \dots \\ 
& \xrightarrow{\hat{G}}\{ f^{p_1,p_2,\dots,p_s}_n(t)\}_{n=0}^\infty,\ p_1,\dots,p_s=1,\dots,k 
\end{split}
\ee  
for every fixed $k,s\in\mbb{N}$ with explicit formulae\footnote{See \cite{6}, (2.19).} for 
\bdis 
f^{p_1,p_2,\dots,p_s}_n(t). 
\edis 

\subsection{} 

In the case of Legendre's orthogonal system 
\be \label{2.3} 
\{ P_n(t)\}_{n=0}^\infty,\ t\in [-1,1]
\ee 
the operator $\hat{G}$ produces, for example, the third generation as follows 
\be \label{2.4} 
\begin{split}
& P_n^{p_1,p_2,p_3}(t)=P_n(u_{p_1}(u_{p_2}(u_{p_3}(t))))\times\prod_{r=0}^{p_1-1}\left|\tilde{Z}(v_{p_1}^r(u_{p_2}(u_{p_3}(t))))\right|\times \\ 
& \prod_{r=0}^{p_2-1}\left|\tilde{Z}(v_{p_2}^r(u_{p_3}(t)))\right|\times 
\prod_{r=0}^{p_3-1}\left|\tilde{Z}(v_{p_3}^r(t))\right|, \\ 
& p_1,p_2,p_3=1, \dots , k,\ t\in [-1,1],\ a=-1,\ l=1, 
\end{split}
\ee 
where 
\be \label{2.5} 
u_{p_i}(t)=\vp_1^{p_i}\left(\frac{\overset{p_i}{\wideparen{T+2}}-\overset{p_i}{T}}{2}(t+1)-\overset{p_i}{T}\right)-T-1,\ i=1,2,3
\ee 
are automorphisms on $[-1,1]$ and 
\be \label{2.6} 
\begin{split}
& v_{p_i}^r(t)=\vp_1^r\left(\frac{\overset{p_i}{\wideparen{T+2}}-\overset{p_i}{T}}{2}(t+1)-\overset{p_i}{T}\right),\ r=0,1,\dots,p_i-1, \\ 
& t\in[-1,1]\Rightarrow u_{p_i}(t)\in [-1,1]\wedge v_{p_i}^r(t)\in[\overset{p_i-r}{T},\overset{p_i-r}{\wideparen{T+2}}]. 
\end{split}
\ee 

\begin{mydef81}
\begin{itemize}
	\item[(a)] Every member of every new $L_2$-orthogonal system 
	\bdis 
	\{P_n^{p_1,p_2,p_3}(t)\}_{n=0}^\infty,\ t\in[-1,1],\ p_1,p_2,p_3=1,\dots,k
	\edis 
	contains the function 
	\bdis 
	\left|\zf\right|_{t=\tau}
	\edis 
	for corresponding $\tau$ since\footnote{See \cite{3}, (9.1), (9.2).} 
	\be \label{2.7} 
	|\tilde{Z}(t)|=\sqrt{\frac{{\rm d}\vp_1(t)}{{\rm d}t}}=\frac{\{1+o(1)\}}{\ln t}\left|\zf\right|,\ t\to\infty; 	
	\ee 
	\item[(b)] Property (a) holds true due to Theorem of the paper \cite{6} for every generation 
	\bdis 
	\{f_n^{p_1,\dots,p_s}(t)\},\ t\in [a,a+2l],\ s\in\mbb{N}. 
	\edis 
\end{itemize}
\end{mydef81} 

\begin{remark}
Our type of \emph{proliferation} of every $L_2$-orthogonal system is in context with the Chumash, Bereishis, 26:12, \emph{Isaac sowed in the land, and in that year reaped a hundredfold, thus had HASHEM blessed him}. 
\end{remark}

\section{Jacob's ladders and following contributions to the theory of Hardy-Littlewood integral}

\subsection{} 

Let us remind that we have introduced the Jacob's ladders 
\be \label{3.1} 
\vp_1(T)=\frac{1}{2}\vp(T) 
\ee  
in our paper \cite{2} (comp. also \cite{3}), where the function $\vp(T)$ is arbitrary solution of the non-linear integral equation\footnote{Also introduced in \cite{2}.} 
\be \label{3.2} 
\int_0^{\mu[x(T)]}\left|\zf\right|^2e^{-\frac{2}{x(T)}t}{\rm d}t=\int_0^T\left|\zf\right|^2{\rm d}t, 
\ee 
where each admissible function $\mu(y)$ generates a solution 
\be \label{3.3} 
y=\vp(T;\mu)=\vp(T),\ \mu(y)\geq 7y\ln y. 
\ee 
We call the function $\vp_1(T)$ Jacob's ladder since analogy with the Jacob's dream in Chumash, Bereishis, 28:12. 

\subsection{} 

Let us remind the Hardy-Littlewood-Ingham formula\footnote{Comp. \cite{7}, (2.1) -- (2.7).} 
\be \label{3.4} 
\int_0^T\left|\zf\right|^2{\rm d}t=T\ln T-(1+\ln2\pi-2c)T+R(T)
\ee 
with the error term $R(T)$ that we write down in the following form 
\be \label{3.5} 
R(T)=\mcal{O}(T^{a+\delta}),\ \frac 14\leq a\leq \frac 13. 
\ee 
In this form the symbol $a$ stands for the least valid exponent of the mentioned type of estimate. 

Next, 83 years after HLI formula (\ref{3.4}), we have proved new result (see \cite{2}). 

\begin{mydef91}
The Hardy-Littlewood integral 
\be \label{3.6} 
J(T)=\int_0^T\left|\zf\right|^2{\rm d}t
\ee 
(see \cite{1}) has, in addition to previously known HLI expression (\ref{3.4}) possessing an unbounded error term at $T\to\infty$, the following infinite number of almost exact representations 
\be \label{3.7} 
\begin{split}
& \int_0^T\left|\zf\right|^2{\rm d}t=\vp_1(T)\ln\{\vp_1(T)\}+ \\ 
& (c-\ln2\pi)\vp_1(T)+c_0+\mcal{O}\left(\frac{\ln T}{T}\right), 
\end{split}
\ee 
where $c_0$ is the constant from the Titschmarsh-Kober-Atkinson formula. 
\end{mydef91} 

\subsection{} 

Further, we have obtained also the following.\footnote{See \cite{2}, (6.2).} 

\begin{mydef92}
\be \label{3.8} 
T-\vp_1(T)\sim(1-c)\pi(T);\ \pi(T)\sim\frac{T}{\ln T},\ T\to\infty, 
\ee  
where the Jacob's ladder can be viewed as the complementary function to the fction $(1-c)\pi(T)$ in the sense 
\be\label{3.9} 
\vp_1(T)+(1-c)\pi(T)\sim T,\ T\to\infty. 
\ee 
\end{mydef92} 

\begin{remark}
The following exchange 
\be \label{3.10}  
\frac{T}{\ln T} \to \pi(T),\ T\to\infty 
\ee  
has been used in the first asymptotic formula (\ref{3.8}).\footnote{See \cite{2}, (6.2).} However, this is completely correct in the asymptotic regions since the prime-number law. At the same time, we can give also a motivation for the exchange (\ref{3.10}). Namely, following 
\be \label{3.11} 
\zeta(s)=\exp\left\{ s\int_2^\infty \frac{\pi(x)}{x(x^s-1)}{\rm d}x\right\},\ \re s>1,  
\ee  
we can get that the function $\pi(x)$ is one element of the set of five functions generating the function $\zeta(s)$. 
\end{remark} 

\subsection{} 

Next, we have proved\footnote{See \cite{7}, (3.4).} the existence of almost linear increments of the Hardy-Littlewood integral (\ref{3.6}). 

\begin{mydef93}
For every fixed natural number $k$ and for every sufficiently big $T>0$ we have 
\be \label{3.12} 
\begin{split}
& \int_{\overset{r-1}{T}}^{\overset{r}{T}}\left|\zf\right|^2{\rm d}t=(1-c)\overset{r-1}{T}+\mcal{O}(T^{a+\delta}), \\ 
& 1-c\approx 0.42,\ r=1,\dots,k,\ T\to\infty, 
\end{split}
\ee  
where, see (\ref{1.5}), (d), 
\be\label{3.13} 
\overset{r}{T}=\overset{r}{T}(T)=\vp_1^{-r}(T). 
\ee 
\end{mydef93} 

\begin{remark}
The existence of linear increments\footnote{Linear with respect to variables $T,\overset{1}{T},\dots,\overset{k}{T}$.} in (\ref{3.12}) is a new phenomenon in the theory of the Riemann's zeta-function. For example, it is true\footnote{See (\ref{3.5}).} that 
\be \label{3.14} 
\int_{\overset{r-1}{T}}^{\overset{r}{T}}\left|\zf\right|^2{\rm d}t=(1-c)\overset{r-1}{T}+\mcal{O}(T^{1/3+\delta}),\ T\to\infty. 
\ee 
\end{remark} 

\begin{remark}
Our formula (\ref{3.12}) represents the result of certain interaction between formulae (\ref{3.4}) and (\ref{3.7}) generated by the Jacob's ladder by means of the property (\ref{1.11}).\footnote{See \cite{7}, (7.5) and (7.6).}
\end{remark}

\subsection{} 

Finally, let us remind our $\zeta$-condition (\ref{1.1}) and, for example, 
\be \label{3.15} 
\lim_{\tau\to\infty}\ln\left\{
\int_{\exp(\frac{x^n+y^n}{z^n}\ln\tau)}^{[\exp(\frac{x^n+y^n}{z^n}\ln\tau)]^1}\left|\zf\right|^2{\rm d}t
\right\}^{\frac{1}{\ln\tau}}\not=1, 
\ee  
see \cite{8}, (6.12), on the class of all Fermat's rationals, which represents two $\zeta$-equivalents of the Fermat-Wiles theorem. 
\begin{remark}
Both mentioned $\zeta$-equivalents represent next contribution to the theory of the Hardy-Littlewood integral (\ref{3.6}). 
\end{remark} 

\section{Formula for certain increments of the Dirichlet's $D(x)$-function} 

\subsection{} 

By makink use of the substitution 
\bdis 
T\to\overset{r-1}{T} 
\edis  
in the Dirichlet's formula (\ref{1.3}) and 
\bdis 
T\to\overset{r}{T} 
\edis  
in our almost exact formula (\ref{3.7}) we obtain the following couple of formulas 
\be \label{4.1} 
D(\overset{r-1}{T})=\overset{r-1}{T}\ln \overset{r-1}{T}+(2c-1)\overset{r-1}{T}+\mcal{O}(T^{1/2}), 
\ee 
\be \label{4.2} 
\begin{split} 
& \int_0^{\overset{r}{T}}\left|\zf\right|^2{\rm d}t=\overset{r-1}{T}\ln \overset{r-1}{T}+(c-\ln2\pi)\overset{r-1}{T}+c_0+\mcal{O}\left(\frac{\ln T}{T}\right), \\ 
& r=1,\dots,k,\ T\to\infty, 
\end{split} 
\ee 
where we have used the properties (\ref{1.11}) and (\ref{1.12}). 

If we subtract (\ref{4.2}) from (\ref{4.1}) we obtain 
\be \label{4.3} 
D(\overset{r-1}{T})-\int_0^{\overset{r}{T}}\left|\zf\right|^2{\rm d}t=(c+\ln2\pi-1)\overset{r-1}{T}+\mcal{O}(T^{1/2}), 
\ee  
and also 
\be \label{4.4} 
D(\overset{r}{T})-\int_0^{\overset{r+1}{T}}\left|\zf\right|^2{\rm d}t=(c+\ln2\pi-1)\overset{r}{T}+\mcal{O}(T^{1/2}), 
\ee 
by the translation $r\to r+1$ in (\ref{4.3}). Now, if we subtract (\ref{4.3}) from (\ref{4.4}), we obtain 
\be\label{4.5} 
\begin{split}
& D(\overset{r}{T})-D(\overset{r-1}{T})=-\int_{\overset{r}{T}}^{\overset{r+1}{T}}\left|\zf\right|^2{\rm d}t= \\ 
& (c+\ln2\pi-1)(\overset{r}{T}-\overset{r-1}{T})+\mcal{O}(T^{1/2}). 
\end{split}
\ee 

\subsection{} 

Next, the almost exact formula (\ref{3.12}) implies 
\be \label{4.6} 
\begin{split}
& \int_{\overset{r}{T}}^{\overset{r+1}{T}}\left|\zf\right|^2{\rm d}t-\int_{\overset{r-1}{T}}^{\overset{r}{T}}\left|\zf\right|^2{\rm d}t= \\ 
& (1-c)(\overset{r}{T}-\overset{r-1}{T})+\mcal{O}(T^{a+\delta}). 
\end{split}
\ee 
Now, it follows from (\ref{4.5}) by (\ref{4.6}) that 
\be \label{4.7} 
\begin{split}
& D(\overset{r}{T})-D(\overset{r-1}{T})=\int_{\overset{r-1}{T}}^{\overset{r}{T}}\left|\zf\right|^2{\rm d}t+ \\ 
& (\ln2\pi)(\overset{r}{T}-\overset{r-1}{T})+\mcal{O}(T^{a+\delta}). 
\end{split} 
\ee  
Finally, by making use the formula (see (\ref{1.12})) 
\be \label{4.8} 
\overset{r}{T}-\overset{r-1}{T}\sim(1-c)\frac{\overset{r}{T}}{\ln \overset{r}{T}}\sim(1-c)\frac{T}{\ln T},\ T\to\infty 
\ee 
in (\ref{4.7}), we obtain the following. 

\begin{mydef51}
\be \label{4.9} 
\begin{split}
& D(\overset{r}{T})-D(\overset{r-1}{T})=\int_{\overset{r-1}{T}}^{\overset{r}{T}}\left|\zf\right|^2{\rm d}t+\mcal{O}\left(\frac{T}{\ln T}\right), \\ 
& r=1,\dots,k,\ k\in\mbb{N}. 
\end{split}
\ee 
\end{mydef51} 

\begin{remark}
Our Remark 1 is based on the estimate 
\bdis 
\mcal{O}\left(\frac{T}{\ln T}\right)
\edis  
of the error term in (\ref{4.9}). 
\end{remark} 

\section{The existence of $D$-equivalents of the Fermat-Wiles theorem}

\subsection{} 

In what follows we shall use, for example, the formula (\ref{4.9}) with $r=1$, that is 
\be \label{5.1} 
D(\overset{1}{T})-D(T)=\int_{T}^{\overset{1}{T}}\left|\zf\right|^2{\rm d}t+\mcal{O}\left(\frac{T}{\ln T}\right),\ T>T_0>0, 
\ee  
where $T_0$ is sufficiently big and 
\be \label{5.2} 
\overset{1}{T}=[T]^1=\vp_1^{-1}(T). 
\ee 
Now, if we put 
\be \label{5.3} 
T=\frac{x}{1-c}\tau,\ \tau\in \left(\frac{1-c}{x}T_0,+\infty\right),\ x>0 
\ee 
into (\ref{5.1}), then we obtain the following statement. 

\begin{mydef52}
\be \label{5.4} 
\begin{split}
& D\left(\left[\frac{x}{1-c}\tau\right]^1\right)-D\left(\frac{x}{1-c}\tau\right)= \\ 
& \int_{\frac{x}{1-c}\tau}^{\left[\frac{x}{1-c}\tau\right]^1}\left|\zf\right|^2{\rm d}t+\mcal{O}\left(\frac{\tau}{\ln \tau}\right) \\ 
& \tau\in (\tau_1(x),+\infty),\ \tau_1(x)=\max\left\{\left(\frac{1-c}{x}\right)^2,(T_0)^2\right\} 
\end{split}
\ee 
for every fixed $x>0$, where, of course, 
\be \label{5.5}  
\frac{1-c}{x}T_0\leq \tau_1(x),\ x>0, 
\ee 
and the constant in the $\mcal{O}$-term depends on $x$. 
\end{mydef52}

Since\footnote{See \cite{8}, (4.6).} 
\be \label{5.6} 
\lim_{\tau\to\infty}\frac{1}{\tau}\int_{\frac{x}{1-c}\tau}^{\left[\frac{x}{1-c}\tau\right]^1}\left|\zf\right|^2{\rm d}t=x 
\ee  
for every fixed $x>0$, then it follows from (\ref{5.4}) the next lemma. 

\begin{mydef53}
\be \label{5.7} 
\lim_{\tau\to\infty}\frac{1}{\tau}\left\{ D\left(\left[\frac{x}{1-c}\tau\right]^1\right)-D\left(\frac{x}{1-c}\tau\right)\right\}=x 
\ee  
for every fixed $x>0$, where (see (\ref{5.2})) 
\be \label{5.8} 
\left[\frac{x}{1-c}\tau\right]^1=\vp_1^{-1}\left(\frac{x}{1-c}\tau\right). 
\ee 
\end{mydef53} 

\subsection{} 

Now, if we use the substitution 
\be \label{5.9} 
x\to \frac{x^n+y^n}{z^n},\ x,y,z\in\mbb{N},\ n\geq 3
\ee 
in (\ref{5.7}), then we obtain the following lemma. 

\begin{mydef54}
\be \label{5.10} 
\lim_{\tau\to\infty}\frac{1}{\tau}\left\{
D\left(\left[\frac{x^n+y^n}{z^n}\frac{\tau}{1-c}\right]^1\right)-D\left(\frac{x^n+y^n}{z^n}\frac{\tau}{1-c}\right)
\right\}=\frac{x^n+y^n}{z^n}
\ee  
for every fixed Fermat's rational 
\bdis 
\frac{x^n+y^n}{z^n}. 
\edis 
\end{mydef54} 

Consequently, we have the following Theorem. 

\begin{mydef11} 
The $D$-condition 
\be \label{5.11} 
\lim_{\tau\to\infty}\frac{1}{\tau}\left\{
D\left(\left[\frac{x^n+y^n}{z^n}\frac{\tau}{1-c}\right]^1\right)-D\left(\frac{x^n+y^n}{z^n}\frac{\tau}{1-c}\right)
\right\}\not=1
\ee 
on the class of Fermat's rationals represents the $D$-equivalent of the Fermat-Wiles theorem. 
\end{mydef11} 

\subsection{} 

Next, let us remind that in the case\footnote{See \cite{8}, (6.6) -- (6.9).} 
\be \label{5.12} 
T=\tau^x,\ x>0,\ \tau\in ((T_0)^{1/x},+\infty)
\ee 
we have obtained the formula 
\be\label{5.13} 
\lim_{\tau\to\infty}\ln\left\{
\int_{\exp(x\ln\tau)}^{[\exp(x\ln\tau)]^1}\left|\zf\right|^2{\rm d}t
\right\}^{\frac{1}{\ln\tau}}=x, 
\ee 
where 
\be \label{5.14} 
\tau\in(\tau_2(x),+\infty),\ \tau_2(x)=\max\{(T_0)^{1/x},T_0\},\ x>0. 
\ee  

Now, if we use the substitution (\ref{5.9}) in (\ref{5.13}), then we obtain by (\ref{5.1}), (\ref{5.12}) -- (\ref{5.14}), the following $D$-equivalent that correspond to the $\zeta$-equivalent in \cite{8}, (6.12). Namely, the following Theorem holds true. 

\begin{mydef12}
The $D$-condition 
\be \label{5.15} 
\lim_{\tau\to\infty}\ln\left\{
D\left(\left[\exp\left(\frac{x^n+y^n}{z^n}\ln\tau\right)\right]^1\right)
-D\left(\exp\left(\frac{x^n+y^n}{z^n}\ln\tau\right)\right)
\right\}\not=1
\ee 
on the class of all Fermat's rationals represents the next $D$-equivalent of the Fermat-Wiles theorem. 
\end{mydef12} 

\section{Further $\zeta$-equivalents of the Fermat-Wiles theorem}

\subsection{} 

We have proved the next formula in our previous paper 
\be \label{6.1} 
\begin{split}
& \int_{\overset{r-1}{T}}^{\overset{r}{T}}\left\{ \sigma(l)\left|\zf\right|^2+(1-c)|S_1(t)|^{2l}\right\}{\rm d}t=(1-c)\sigma(l)\overset{r}{T}+ \\ 
& \mcal{O}\left(\frac{T}{\ln^2 T}\right),\ r=1,\dots,k,\ k,l\in\mbb{N} 
\end{split}
\ee 
for every fixed $k$ and $l$ (see\cite{8}, (3.19)). Here we shall use the formula (\ref{6.1}), $r=2$, i. e. 
\be \label{6.2} 
\begin{split} 
& \int_{\overset{1}{T}}^{\overset{2}{T}}\left\{ d(l)\left|\zf\right|^2+(1-c)|S_1(t)|^{2l}\right\}{\rm d}t=(1-c)d(l)\overset{2}{T}+ \\ 
& \mcal{O}\left(\frac{T}{\ln^2 T}\right) 
\end{split} 
\ee 
for example, where 
\be \label{6.3} 
\begin{split}
& \overset{1}{T}=\overset{1}{T}(T)=\vp_1^{-1}(T), \\ 
& \overset{2}{T}=\overset{2}{T}(T)=\vp_1^{-2}(T). 
\end{split}
\ee 
Now we put 
\be \label{6.4} 
T=\vp_1^2(\tau);\ \tau=\vp_1^{-2}(T)=\overset{2}{T} 
\ee 
in (\ref{6.3}) to obtain 
\be \label{6.5} 
\overset{1}{T}=\vp_1^{-1}(T)=\vp_1(\tau). 
\ee  
Asymptotic formulae in (\ref{1.12}) imply by (\ref{6.4}) and (\ref{6.5}) the following 
\be \label{6.6} 
\tau\sim\vp_1(\tau)\sim\vp_1^2(\tau),\ \tau\to\infty. 
\ee 
Next, by (\ref{6.6}) we have (see (\ref{6.2})) 
\be \label{6.7} 
\left. \mcal{O}\left(\frac{T}{\ln^2 T}\right)\right|_{T=\vp_1^2(\tau)}=\mcal{O}\left(\frac{\vp_1^2(\tau)}{[\ln\vp_1^2(\tau)]^2}\right)=\mcal{O}\left(\frac{\tau}{\ln^2 \tau}\right). 
\ee 
Consequently we obtain (see (\ref{6.5}) -- (\ref{6.7})) the following lemma. 

\begin{mydef55}
\be \label{6.8} 
\begin{split}
& \int_{\vp_1(\tau)}^\tau\left\{
\sigma(l)\left|\zf\right|^2+(1-c)|S_1(t)|^{2l}
\right\}{\rm d}t= \\ 
& (1-c)\sigma(l)\tau+\mcal{O}\left(\frac{\tau}{\ln^2 \tau}\right),\ \tau\in (\overset{2}{T}(T_0),+\infty). 
\end{split}
\ee 
\end{mydef55} 

\subsection{} 

Now, if we put 
\be \label{6.9} 
\tau=\frac{x}{1-c}\rho,\ \rho\in \left(\frac{(1-c)d}{x}\overset{2}{T}(T_0),+\infty\right),\ x>0
\ee 
into (\ref{6.8}) then we obtain the following lemma. 

\begin{mydef56}
\be \label{6.10} 
\begin{split}
& \frac{1}{\rho}\int_{\vp_1(\frac{x}{(1-c)\sigma(l)}\rho)}^{\frac{x}{(1-c)\sigma(l)}\rho}\left\{
\sigma(l)\left|\zf\right|^2+(1-c)|S_1(t)|^{2l}
\right\}{\rm d}t= \\ 
& x+\mcal{O}\left(\frac{1}{\ln^2\rho}\right),\ x>0, \\ 
& \rho\in(\rho_1,+\infty),\ \rho_1=\max\left\{\left(\frac{(1-c)d}{x}\right)^2,(\overset{2}{T}(T_0))^2\right\}. 
\end{split}
\ee 
\end{mydef56} 

Next, if we use the substitution 
\bdis 
x\mapsto \frac{x^n+y^n}{z^n} 
\edis 
in (\ref{6.10}) then we obtain the following statement. 

\begin{mydef57}
We have for every fixed Fermat's rational $\frac{x^n+y^n}{z^n}$ that 
\be \label{6.11} 
\begin{split}
& \lim_{\rho\to\infty}\frac{1}{\rho}\int_{\vp_1(\frac{x^n+y^n}{z^n}\frac{\rho}{(1-c)\sigma(l)})}^{\frac{x^n+y^n}{z^n}\frac{\rho}{(1-c)\sigma(l)}}\left\{
\sigma(l)\left|\zf\right|^2+(1-c)|S_1(t)|^{2l}
\right\}{\rm d}t= \\ 
& \frac{x^n+y^n}{z^n}. 
\end{split}
\ee 
\end{mydef57}

Now the theorem follows. 

\begin{mydef13}
The $\zeta$-condition 
\be \label{6.12} 
\begin{split}
& \lim_{\rho\to\infty}\frac{1}{\rho}\int_{\vp_1(\frac{x^n+y^n}{z^n}\frac{\rho}{(1-c)\sigma(l)})}^{\frac{x^n+y^n}{z^n}\frac{\rho}{(1-c)\sigma(l)}}\left\{
\sigma(l)\left|\zf\right|^2+(1-c)|S_1(t)|^{2l}
\right\}{\rm d}t\not=1 
\end{split}
\ee 
on the class of all Fermat's rationals represents the next $\zeta$-equivalent of the Fermat-Wiles theorem, where 
\bdis 
S_1(t)=\frac{1}{\pi}\int_0^t\arg\zf{\rm d}t 
\edis  
for completeness. 
\end{mydef13} 

\section{Concluding remarks} 

\subsection{} 

Let 
\be \label{7.1} 
\begin{split}
& L_1=\lim_{\tau\to\infty}\frac{1}{\tau}\int_{a\frac{\tau}{1-c}}^{[a\frac{\tau}{1-c}]^1}\left|\zf\right|^2{\rm d}t,\ a>0, \\ 
& L_2=\lim_{\tau\to\infty}\frac{1}{\tau}\int_{\frac{x^n+y^n}{z^n}\frac{\tau}{1-c}}^{[\frac{x^n+y^n}{z^n}\frac{\tau}{1-c}]^1}\left|\zf\right|^2{\rm d}t, \\ 
& L_3=\lim_{\tau\to\infty}\frac{1}{\tau}\int_{a\frac{x^n+y^n}{z^n}\frac{\tau}{1-c}}^{[a\frac{x^n+y^n}{z^n}\frac{\tau}{1-c}]^1}\left|\zf\right|^2{\rm d}t. 
\end{split}
\ee  
Of course, it is true that 
\be \label{7.2} 
L_1L_2=L_3, 
\ee  
(see \cite{8}, (4.11)). 

Now, we have the following: 
\begin{itemize} 
\item[(A)] $L_2\not=1 \ \Rightarrow \ aL_2\not=a \ \Rightarrow \ L_3\not=a,\ \forall a>0$; 
\item[(B)] $L_3\not=a \ \Rightarrow \ L_1L_2\not=a \ \Rightarrow \ aL_2\not=a \ \Rightarrow \ L_2\not=1,\ \forall a>0$; 
\end{itemize}  
see (\ref{7.2}). Consequently, it is true the following. 

\begin{remark}
\be \label{7.3} 
L_2\not=1 \ \Leftrightarrow \  L_3\not=a
\ee  
for every fixed $a>0$. 
\end{remark}

\subsection{} 

Next, we have by Lemma 3, see (\ref{5.7}), the following. 

\begin{remark}
\begin{itemize}
	\item[(A)] In the case 
	\be \label{7.4} 
	x\to D(x) 
	\ee 
	we obtain 
	\be \label{7.5} 
	\lim_{\tau\to\infty}\frac{1}{\tau}\left\{
	D\left(\left[D(x)\frac{\tau}{1-c}\right]^1\right)-D\left(D(x)\frac{\tau}{1-c}\right)
	\right\}=D(x), 
	\ee 
	that gives us the return to the Dirichlet's function $D(x)$; 
	\item[(B)] in the case of Euler's Gamma function 
	\be \label{7.6} 
	x\to\Gamma(x),\ \Gamma(\Gamma(x)), \dots , x>0 
	\ee 
	we obtain, for example,  
	\be \label{7.7} 
	\begin{split}
	& \lim_{\tau\to\infty}\frac{1}{\tau}\left\{
	D\left(\left[\Gamma(\Gamma(x))\frac{\tau}{1-c}\right]^1\right)-
	D\left(\Gamma(\Gamma(x))\frac{\tau}{1-c}\right)=\Gamma(\Gamma(x))
	\right\}, 
	\end{split}
	\ee  
	that implies a kind of continuation of the formula (\ref{5.7}) into infinite set of other possibilities. 
\end{itemize}
\end{remark}

I would like to thank Michal Demetrian for his moral support of my study of Jacob's ladders.

\end{document}